\documentclass[reqno,11pt]{amsart}

\usepackage{amssymb,amsthm}
\usepackage{amsmath}
\usepackage{amsfonts}
\usepackage{dsfont}

\usepackage[a4paper,  margin=3.7cm]{geometry}

\usepackage[active]{srcltx}
\makeatletter\@addtoreset{equation}{section}\makeatother

\newtheorem{theorem}{Theorem}[section]

\newtheorem{remark}[theorem]{Remark}
\numberwithin{equation}{section}

\title{Self-Regulation in Infinite Populations with Fission-Death Dynamics}

\author{ Yuri  Kozitsky}

\address{Instytut Matematyki, Uniwersytet Marii Curie-Sk{\l}odowskiej, 20-031 Lublin, Poland}
\email{jkozi@hektor.umcs.lublin.pl}

\author{Agnieszka Tana{\'s}}

\address{Politechnika Lubelska, 20-618 Lublin, Poland}

\begin{document}


\begin{abstract}

The evolution of an infinite population of interacting point
entities placed in $\mathbb{R}^d$ is studied. The elementary
evolutionary acts are death of an entity with rate that includes a
competition term and independent fission into two entities. The
population states are probability measures on the corresponding
configuration space and the result is the construction of the
evolution of states in the class of sub-Poissonian measures, that
corresponds to the lack of clusters in such states. This is
considered as a self-regulation in the population due to
competition.

\end{abstract}

\maketitle

\section{Introduction}
\label{S1}
\subsection{Regulating population dynamics}
Simple population dynamics models are mostly based on two
evolutionary acts: disappearance (death) of an entity and
procreation, in the course of which new entities join the
population. A commonly accepted viewpoint on the evolution of a
finite population of this kind is that it either dies out or grows
ad infinitum being unhampered by habitat restrictions. Clearly, such
restrictions can only be ignored if the population size is small,
i.e., at the early stage of its development. In developed
populations, environmental restrictions force the entities to
compete with each other -- a crowding effect. In the mentioned
models, this effect manifests itself in a state-dependent increment
of the death toll. In Verhulst's phenomenological theory based on
the equation $\frac{d}{dt} N = \lambda N - (\mu + \alpha N)N$, such
an increment is $\alpha N$. Here $N=N(t)$ is the (expected) number
of entities at time $t$, and positive $\lambda$ and $\mu$ are the
intrinsic procreation and death rates, respectively. Later on, Pearl
and Reed rewrote this in the form of the logistic growth equation
$\frac{d}{dt} N = r N ( 1 - N/K)$ with $r= \lambda - \mu$ and $K= r/
\alpha$. The latter parameter gives rise to the notion of
\emph{carrying capacity} as the solution $N\equiv K$ is a stationary
one, to which $N(t)$ tends in the limit $t\to +\infty$. Since then,
this notion is used in the theory of biological populations, see
Introduction in \cite{HJK}, and not only in the context of the
competition caused crowding effect. For instance, in the
Galton-Watson model with binary fission considered in \cite{Fima},
the probability of fission of a member of generation $n$ consisting
of $Z_n$ entities was taken to be $K/(K+Z_n)$. Thereby, the
constructed process gets super-- or subcritical under or over the
level $K$, respectively. This aspect of the theory may be viewed as
a phenomenological (mean-field-like) way of regulating the
population dynamics. Here {\it regulating} means preventing the
population from infinite growth and {\it mean-field} corresponds to
imitating interactions as state-dependent external actions (fields),
cf. \cite[Sect. 13]{Simon}.

In the theory of populations with interactions explicitly taken into
account, a usual assumption is that each entity interacts mostly (or
even entirely) with the subpopulation located in a compact subset of
the habitat. Then the local structure of the population is
determined by the network of such interactions. Since a finite
population occupies a compact set, it is always local as each of its
members has a compact neighborhood containing the whole remaining
population. Thus, in order to understand the global behavior of
populations of this kind, one should take them infinite. In the
statistical mechanics of interacting physical particles developed
from phenomenological thermodynamics, this conclusion had led to the
concept of the infinite-volume limit, see, e.g., \cite[pp.
5,6]{Simon}. In this note, and in the accompanying paper \cite{KT}
where all the technical details are presented, we introduce an
individual-based model of an infinite population of point entities
placed in $\mathds{R}^d$ which undergo binary fission and death
caused also by crowding (local competition). Its aim is to to
demonstrate that the local competition -- interaction explicitly
taken into account -- can produce a global regulating effect. Here,
however, one has to make precise the very notion of regulation as
the considered population is already infinite. Instead of
characterizing it by the number of constituents, we will look at the
spatial distribution of the population by comparing it with the
distribution governed by a Poisson law.

\subsection{Presenting the result}
Similarly as in \cite{KK}, we  deal with the phase space $\Gamma$
consisting of all locally finite subsets  $\gamma
\subset\mathds{R}^d$, called {\it configurations}. Local finiteness
means that $\gamma_\Lambda := \gamma\cap \Lambda$ is finite whenever
$\Lambda \subset\mathds{R}^d$ is compact. For compact $\Lambda$ and
$n\in \mathds{N}_0$, we then set $\Gamma^{\Lambda,n} = \{\gamma\in
\Gamma: |\gamma_\Lambda|=n\}$, where $|\cdot|$ denotes cardinality,
and equip $\Gamma$ with the $\sigma$-field $\mathcal{B}(\Gamma)$
generated by all such $\Gamma^{\Lambda,n}$. This allows one to
consider probability measures on $\Gamma$ as states of the system.
In a Poisson state, the entities are independently distributed over
$\mathds{R}^d$. A homogeneous Poisson measure $\pi_\varkappa$ with
intensity $\varkappa>0$ is characterized by its values on
$\Gamma^{\Lambda,n}$ given by the following expression
\begin{equation}
  \label{1}
  \pi_\varkappa (\Gamma^{\Lambda,n}) = \frac{(\varkappa
  |\Lambda|)^n}{n!} \exp\left(- \varkappa |\Lambda| \right),
\end{equation}
where $|\Lambda|$ stands for the Lebesgue measure of $\Lambda$. Note
that $\pi_\varkappa (\Gamma_0)=0$, for all $\varkappa>0$, where
$\Gamma_0$ is the set of all finite configurations. Let
$\mathcal{P}(\Gamma)$ be the set of all probability measures on
$\Gamma$. We say that a given $\mu\in \mathcal{P}(\Gamma)$ is {\it
sub-Poissonian} if, for each compact $\Lambda$, all $n\in
\mathds{N}_0$ and some $\varkappa>0$, the following holds
\begin{equation}
  \label{2}
  \mu (\Gamma^{\Lambda,n}) \leq \pi_\varkappa(\Gamma^{\Lambda,n}).
\end{equation}
It is believed that sub-Poissonian states are characterized by the
lack of {\it clustering}, typical to procreating populations with
noninteracting (noncompeting) constituents, see the corresponding
discussion in \cite{KK}.

In dealing with states on $\Gamma$, one employs {\it observables} --
appropriate functions $F:\Gamma \to \mathds{R}$. Their evolution is
obtained by solving the Kolmogorov equation
\begin{equation*}
\frac{d}{dt} F_t = L F_t, \qquad F_{t}|_{t=0} = F_0, \qquad t>0,
\end{equation*}
in which the operator $L$ specifies the model. The model which we
introduce here is based on the following evolutionary acts: (a) an
entity located at $x$ dies with rate (probability per unit time)
$m(x) + \sum_{y\in \gamma \setminus x} a(x-y)$, where $m(x)\geq 0$
corresponds to a per se mortality and $a\geq0$ is the competition
kernel; (b) an entity located at $x$ undergoes fission, with two
offsprings going to $y_1$ and $y_2$ with rate $b(x|y_1,y_2)$.
According to this, the operator $L$ takes the form
\begin{gather}
  \label{2a}
 (LF)(\gamma) = \sum_{x\in \gamma}\left( m (x) + \sum_{y\in \gamma\setminus x} a (x-y)\right)
\left[ F(\gamma\setminus x) - F(\gamma)
 \right] \\[.2cm] \nonumber + \sum_{x\in \gamma}
 \int_{(\mathds{R}^d)^2} b(x|y_1 , y_2)\left[ F(\gamma\setminus x \cup\{y_1, y_2\}) - F(\gamma)
 \right] dy_1 d y_2.
\end{gather}
In expressions like $\gamma \cup x$, we treat $x$ as the singleton
$\{x\}$. Note that the fission rate is state-independent. Regarding
$a$, $m$ and $b$ we assume that: (a) $a:\mathds{R}^d \to
[0,+\infty)$ is a piece-wise continuous function such that $a(x)=0$
whenever $|x|>r$ for some positive $r<\infty$; (b) $m:\mathds{R}^d
\to [0,+\infty)$ is measurable and bounded; (c) the fission kernel
is translation invariant in the sense that $b(x+z|y_1 +z, y_2 + z) =
b(x|y_1, y_2)$ holding for all $z\in \mathds{R}^d$; (d) the function
$\beta:\mathds{R}^d \to [0,+\infty)$ defined by
\begin{equation}
  \label{3A}
  \beta (y_1 - y_2) = \int_{\mathds{R}^d} b(x|y_1 , y_2) d x
\end{equation}
is piece-wise continuous and such that $\beta(x)=0$ whenever $|x|>R$
for some positive $R<\infty$; (e) $\beta (x) = \beta (-x)$ for all
$x\in \mathds{R}^d$ and the following holds
\begin{equation}
  \label{3A}
\int_{\mathds{R}^d}  \beta (y)dy = \int_{(\mathds{R}^d)^2} b(x|y_1 ,
y_2) d y_1 d y_2 =:\langle b \rangle<\infty.
\end{equation}
Note that the translation invariance and the finite-range property
are imposed here only to make the presentation of the model and the
results as simple as possible. The version studied in \cite{KT} is
characterized by less restrictive conditions. Note also that we do
not exclude the case where $b$ is a distribution. For instance, by
setting
$$b(x|y_1, y_2)= \frac{1}{2} \left(\delta (x-y_1) + \delta (x-y_2)
\right)\beta (y_1-y_2),$$ we obtain the Bolker-Pacala model
\cite{KK} as a particular case of our model.
\begin{remark}
  \label{1rk}
The function $\beta$ describes the dispersal of siblings, which
compete with each other. As in the Bolker-Pacala model, here the
following situations may occur:
\begin{itemize}
  \item \emph{Short dispersal:} there exists $\omega >0$
  such that $a(x)
   \geq \omega \beta(x)$ for all $x\in \mathds{R}^d$; corresponds to
   $R\leq r$.
  \item \emph{Long dispersal:} for each $\omega >0$, there
  exists $x\in \mathds{R}^d$ such that $ a(x)
  < \omega \beta(x)$; corresponds to
   $R> r$.
\end{itemize}
\end{remark}
The direct use of $L$ as a linear operator in an appropriate Banach
space is possible only if one restricts the consideration to states
on $\Gamma_0$, see \cite[Sect. 3]{KT}. Otherwise, the sums in
(\ref{2a}) -- taken over infinite configurations -- may not exist.
In view of this, we proceed as follows. Let $C_0 (\mathds{R}^d)$
stand for the set of all continuous real-valued functions with
compact support. Then the map
\begin{eqnarray*}
& & \Gamma \ni \gamma \mapsto F^\theta (\gamma)  :=  \prod_{x\in
\gamma} (1+ \theta (x)), \qquad \theta \in \varTheta, \\[.2cm]
& & \varTheta := \{\theta\in C_0 (\mathds{R}^d): \theta(x) \in
(-1,0]\},
\end{eqnarray*}
is measurable and satisfies $0<F^\theta (\gamma) \leq 1$ for all
$\gamma$. It is possible to show, see \cite{KT}, that $\{F_\theta:
\theta \in \varTheta\}$ is a measure defining class. That is, for
each two $\mu,\nu \in \mathcal{P}(\Gamma)$, it follows that
$\mu=\nu$ whenever $\mu(F^\theta):=\int F^\theta d \mu = \int
F^\theta d \nu=:\nu (F^\theta)$ holding for all such $F^\theta$.
Moreover, under the mentioned above assumption that both $a$ and $b$
in (\ref{2a}) have finite range, $(L F^\theta)(\gamma)$ can be
calculated for each $\gamma\in \Gamma$ and $\theta \in \varTheta$.
We prove that, for each $\mu_0 \in \mathcal{P}_{\rm exp}(\Gamma)$,
there exists the map $[0,+\infty) \ni t \mapsto
\mu_t\in\mathcal{P}_{\rm exp}(\Gamma)$ such that
$\mu_t|_{t=0}=\mu_0$, the map $(0,+\infty) \ni t \mapsto \mu_t (
F^\theta)$ is continuously differentiable  and the following holds
\begin{equation}
  \label{4a}
  \frac{d}{dt} \mu_t (F^\theta) = \mu_t (L F^\theta).
 \end{equation}
Here $\mathcal{P}_{\rm exp}(\Gamma)$ is a class of measures each
element of which is sub-Poissonian, see below, and such that $\mu (L
F^\theta)<\infty$.

\section{The Result}

For the Poisson measure as in (\ref{1}), it follows that
$$\pi_\varkappa (F^\theta) = \exp\left(\varkappa \int_{\mathds{R}^d} \theta (x) d x
\right).$$ Having this in mind we introduce the class of measures
$\mathcal{P}_{\rm exp}(\Gamma)$ by the condition that, for each
$\mu\in\mathcal{P}_{\rm exp}(\Gamma)$, $\mu(F^\theta)$ can be
continued to an exponential type entire function of $\theta \in
L^1(\mathds{R}^d)$. It can be shown that $\mu\in\mathcal{P}_{\rm
exp}(\Gamma)$ if and only if $\mu (F^\theta)$ ƒis written in the
form
\begin{equation}
 \label{6b}
 \mu(F^\theta) = 1 + \sum_{n=1}^\infty \frac{1}{n!} \int_{(\mathds{R}^d)^n} k^{(n)}_\mu (x_1 , \dots , x_n) \theta (x_1) \cdots \theta(x_n) d x_1 \cdots d x_n,
\end{equation}
where $k^{(n)}_\mu$ is the $n$-th order \emph{correlation function}
of $\mu$. Each $k^{(n)}_\mu$ is a symmetric positive element of
$L^\infty((\mathds{R}^d)^n)$ satisfying the \emph{Ruelle bound}, cf.
\cite{Ruelle},
\begin{equation}
 \label{6c}
k^{(n)}_\mu (x_1, \dots , x_n) \leq  \varkappa^{ n}, \qquad n \in
\mathds{N},
\end{equation}
holding with some $\varkappa >0$. Note that (\ref{6c}) readily
yields (\ref{2}). Note also that states of thermal equilibrium of
systems of interacting physical particles satisfy (\ref{6c}), see
\cite{Ruelle}. By means of $k^{(n)}_\mu$ one can define the function
$k_\mu : \Gamma_0 \to \mathds{R}$ by setting $k_\mu(\{x_1 , \dots,
x_n\}) = k^{(n)}_\mu(x_1 , \dots, x_n)$, $n\in \mathds{N}$. Let us
consider the Banach space $\mathcal{K}_\alpha$ of such functions
equipped with the norm $$\|k\|_\alpha = \sup_{n\geq 0}
\|k^{(n)}\|_{L^\infty((\mathds{R}^d)^n)} \exp(-\alpha n), \quad
\alpha \in \mathds{R}, $$ and with the usual point-wise linear
operations. Clearly, $\|k\|_{\alpha'} \leq \|k\|_{\alpha}$ whenever
$\alpha' > \alpha$, which yields that
\begin{equation}
  \label{6CC}
 \mathcal{K}_\alpha \hookrightarrow  \mathcal{K}_{\alpha'}, \qquad
 \alpha <
 \alpha'.
\end{equation}
Hence, $\{\mathcal{K}_\alpha\}_{\alpha \in \mathds{R}}$ form an
ascending scale of Banach spaces. In each $\mathcal{K}_\alpha$, one
defines the unbounded linear operator $(L^\Delta,
\mathcal{D}_\alpha)$ by setting $\mathcal{D}_\alpha =\{ k\in
\mathcal{K}_\alpha : L^\Delta k \in \mathcal{K}_\alpha\}$, where the
action of $L^\Delta$ on $k_\mu$ is calculated from the formula, cf.
(\ref{6b}),
\begin{eqnarray}
 \label{6D}
 \mu(L F^\theta) & = & 1 \\[.2cm] \nonumber & + &
 \sum_{n=1}^\infty \frac{1}{n!} \int_{(\mathds{R}^d)^n} (L^\Delta k_\mu)^{(n)} (x_1 , \dots , x_n) \theta (x_1) \cdots \theta(x_n) d x_1 \cdots d x_n,
\end{eqnarray}
Then with the help of (\ref{6b}) and (\ref{6D}) the evolution $\mu_0
\to \mu_t$ is obtained by employing the correlation functions in the
following three steps:
\begin{itemize}
  \item[(a)] Constructing $k_0\to k_t$ for $t< T<\infty$ by solving the
corresponding evolution equation
\begin{equation}
  \label{6CD}
\frac{d}{dt} k_t = L^\Delta k_t, \qquad k_t|_{t=0} = k_{\mu_0}.
\end{equation}
\item[(b)] Proving that $k_t$ is the correlation function of a unique
$\mu_t \in \mathcal{P}_{\rm exp}(\Gamma)$. \item[(c)] Continuing
$k_t$ to all $t>0$.
\end{itemize}
To perform step (a), for each $\alpha_0 \in \mathds{R}$ and
$\alpha_1
> \alpha_0$, we construct a family of operators $Q_{\alpha_1 \alpha_0} (t)$,
$t\in [0, T(\alpha_1, \alpha_0))$. Here $T(\alpha_1, \alpha_0) =
(\alpha_1 - \alpha_0)/ \tau(\alpha_1)$ with certain (dependent on
the model parameters and explicitly found) function $\tau(\alpha)$.
Each $Q_{\alpha_1 \alpha_0} (t)$ acts as a bounded operator from
$\mathcal{K}_{\alpha_0}$ to $\mathcal{K}_{\alpha_1}$, cf.
(\ref{6CC}). Then the (classical) solution of (\ref{6CD}) with $k_0
\in \mathcal{K}_{\alpha_0}$ is obtained in the form $k_t =
Q_{\alpha_1 \alpha_0} (t) k_0$, $t < T(\alpha_1, \alpha_0)$. The
important peculiarities of this solution are: (i) the function
$\tau(\alpha)$ is rapidly increasing, which means that the time
interval shrinks to zero as $\alpha \to +\infty$;  (ii) as $t$
increases, $k_t$ passes to an ever-larger space, cf. (\ref{6CC});
(iii) the solution $k_t$ need not be a correlation function of any
state. In view of (i) and (ii), the direct continuation of $k_t$ to
all $t>0$ is impossible.

To perform step (b) we use a special cone $\mathcal{K}^\star_\alpha
\subset \mathcal{K}_\alpha$ (explicitly constructed, see eq. (4.11)
in \cite{KT}) such that $k\in \mathcal{K}_\alpha$ is the correlation
function of a unique $\mu\in \mathcal{P}_{\rm exp} (\Gamma)$ if and
only if $k\in \mathcal{K}^\star_\alpha$. Then we prove that the
solution mentioned above lies in $\mathcal{K}^\star_{\alpha_1}$ for
all $t < T(\alpha_1, \alpha_0)/3$. Along with the identification of
$k_t$ as a correlation function, this yields also that $k_t\in
\mathcal{K}_{\alpha_t}$ with $\alpha_t = \alpha_0 + c t$. Here
$\alpha_0$ is chosen to be such that $k_0 \in
\mathcal{K}_{\alpha_0}$ and $\alpha_0>-\log \omega$ with $\omega$ as
in Remark \ref{1rk}. One can take $c=0$ if $m_*:=\inf_{x\in
\mathds{R}^d}m(x)> \langle b \rangle$, where the latter is the same
as in (\ref{3A}). In the short dispersal case, one can take $c=0$
already for $m_* = \langle b \rangle$. For $c=0$, the solution stays
in the same space and hence can be continued to all $t>0$ by
repeating the above construction. This is not the case if $c>0$.
Then the solution passes to an ever-larger space, but with a much
slower increase than in the construction made in step (a). This
allows one to prove that $k_t \in \mathcal{K}_{\alpha_t}$ for all
$t>0$ also for positive $c$, and hence to perform step (c). Note
that in the essentially different cases of short and long dispersal
the qualitative difference of the corresponding dynamics appear only
at the borderline case of $m_*=\langle b \rangle$. This may mean
that the dispersal range affects finer properties of the
corresponding system.

As the result, under the assumptions on $a$, $m$ and $b$ made above
we prove the following statement, see \cite[Theorem 4.1 and
Corollary 4.2]{KT}.
\begin{theorem}
  \label{1tm}
There exist $c\in \mathds{R}$ and $\omega >0$ such that, for each
$\mu_0\in \mathcal{P}_{\rm exp}(\Gamma_0)$, there exists a unique
map $[0,+\infty) \ni t \mapsto k_t \in
\mathcal{K}^{\star}_{\alpha_t}$ with $\alpha_t = \alpha_0 + ct$ and
$\alpha_0> - \log \omega$ such that $k_0=k_{\mu_0}\in
\mathcal{K}^\star_{\alpha_0}$, which has the following properties:
\begin{itemize}
  \item[(i)]
For each $T>0$ and all $t\in [0,T)$, the map $$ [0,T)\ni t \mapsto
k_t \in \mathcal{K}_{\alpha_t}  \subset
 \mathcal{D}_{\alpha_T} \subset \mathcal{K}_{\alpha_T}$$
is continuous on $[0,T)$ and continuously differentiable on $(0,T)$
in $\mathcal{K}_{\alpha_T}$.
\item[(ii)] For all $t\in (0,T)$, it satisfies
 $\frac{d}{dt}k_t = L^\Delta k_t$.
\end{itemize}
\end{theorem}
By Theorem \ref{1tm} the evolution $\mu_0\to \mu_t$ in question is
obtained by identifying $\mu_t$ by its values on $F^\theta$ with the
help of (\ref{6b}) and by the evolution $k_0\to k_t$ constructed
therein. Then the validity of (\ref{4a}) follows by (\ref{6D}).
Since $\mu_t\in \mathcal{P}_{\rm exp}(\Gamma)$ for all $t>0$, this
evolution preserves the sub-Poissonicity of the states and hence the
self-regulation -- in the above-mentioned sense -- takes place. Like
in the Bolker-Pacala model, see the corresponding discussion in
\cite{KK}, in our case it can be shown that the self-regulation of
this kind does not hold for $a\equiv 0$.

\section*{Acknowledgment}
The authors are grateful to Krzysztof Pilorz for valuable
discussions. Yuri Kozitsky was supported by National Science Centre,
Poland, grant 2017/25/B/ST1/00051 that is cordially acknowledged.





\end{document}